\documentclass[12pt]{article}
\usepackage{amsmath, amsthm, amscd, amsfonts, amssymb, graphicx, color}
\usepackage[bookmarksnumbered, colorlinks, plainpages]{hyperref}
\newtheorem{theorem}{Theorem}[section]

\newtheorem{proposition}[theorem]{Proposition}
\newtheorem{corollary}{Corollary}[theorem]
\theoremstyle{definition}
\newtheorem{definition}[theorem]{Definition}

\theoremstyle{remark}
\newtheorem{remark}{Remark}[section]
\numberwithin{equation}{section}

\begin{document}
	
	\setcounter{page}{1}
	
	
	\begin{center}
		{\Large \textbf{On Some Characterizations of General $s$-Convex Functions}}
		
		\bigskip

		\textbf{Musavvir Ali$^{a,*}$} and \textbf{Ehtesham Akhter$^b$}\\
		\textbf{}  \textbf{}\\

		{\small $^{a,~b}$  Department of Mathematics,\\ Aligarh Muslim University, Aligarh-202002, India}
	\end{center}
	\noindent
	{\footnote\\
		E-mail addresses: musavvir.alig@gmail.com (M. Ali),\\ ehteshamakhter111@gmail.com (E. Akhter).}
	\bigskip

	{\abstract
		\noindent
		It is established that general $s$-convex functions are a new class of generalized convex functions. In a similar vein, a new class of general $s$-convex sets is introduced, which are generalizations of $s$-convex sets. Additionally, certain fundamental characteristics of general $s$-convex functions are discussed for both general cases and differentiable situations. Aside from that, the general $s$-convexity is used to define and demonstrate the sufficient criteria for optimality for both unconstrained and inequality-constrained programming.}\\
	
	
	\section{Introduction} Due to the significance of convexity and generalized convexity in the investigation of optimality to resolve mathematical programming, researchers focused heavily on the generalized convex functions. For instance, Earlier works by H. Hudzik et al. \cite{Hudzik} presented two types of $s$-convexity \{$s \in (0,1)$\} and demonstrated that Whenever \{$s \in (0,1)$\}, the second notion is fundamentally stronger than the $s$-convexity in the first sense. E-convex sets and $E$-convex functions are a class of sets and a class of functions introduced by Youness \cite{Youness}. Functions by loosening up how convex sets and convex functions are defined. X.M. Yang \cite{X} offered a few instances for E.A. paper \cite{Youness} by Youness and improved it. More findings about generalized E-convex functions, location cite \cite{Duca,Fulga} and any references within that are directly related.
	
	These kinds of generalized convex functions have recently attracted a lot of attention from researchers. Particularly for	study of the $b$-invex function. Consider X.J. Long et al. \cite{Long} talked about a group of functions. A generalization of semi-preinvex functions and $b$-vex functions is known as semi-$b$-preinvex functions. A class of functions known as $E$-$b$-vex functions, which is defined as an extension of $b$-vex functions and $E$-vex functions, was introduced by Yu-Ru Syau et al. \cite{Syau}. A new family of functions known as approximately $b$-invex functions was studied by Emam \cite{Emam}, some of its characteristics were explored, and suitable optimality requirements for nonlinear programming utilising these functions were obtained. A family of sub-$b$-convex sets and novel generalized sub-b-convex functions were examined by M.T. Chao et al. in 2012 \cite{Chao}, who also provided the necessary conditions for the optimality of both unconstrained and inequality-constrained sub-$b$-convex programming. For More details on generalized convex functions can be found in \cite{Mishra,Z}. Through their research, generalized convex functions like the $b$-invex function were developed. Through their research, generalized convex functions like the $b$-invex function were developed. In the meantime, we now discover a class of generalized convex functions that is more comprehensive than these two types of generalized convex functions and which is not sub-$b$-convex but which shares some of their characteristics. Because sub-$b$-convexity and $s$-convexity are extensions of convexity, they pique our curiosity in fresh study, thus we focus on it.
	
	This paper's goal is to introduce a new class of generalized convex functions known as general $s$-convex functions and explore some characteristics of the class of functions meeting general $s$-convexity. It is motivated by research works \cite{Hudzik,Chao,Stancu}. Additionally, we offer adequate guidelines for both of their optimality programming that is obtained under general $s$-convexity is both unrestricted and inequality-restricted.
	
	The rest of this essay is divided into the following sections. Section 2, which further expands on the idea of general $s$-convexity, introduces a new class of functions known as general $s$-convex functions. A new class of sets termed general $s$-convex sets is consequently introduced, as well as some characteristics of the general $s$-convex function
	developed, including general $s$-convex sets. In Section 3, we present fresh general $s$-convex programming and establish sufficient criteria for optimality within the general $s$-convexity. 
	
	\section{Basic Results}
	We reviewed the definitions of sub-$b$-convexity, $s$-convexity and sub-$b$-$s$-convexity of function at the beginning of this section. According to M.T. Chao et al. \cite{Chao}, the class of sub-$b$-convex functions is as follows. Let $\aleph$ be a convex set in $\mathbf{R^m}$ which is non-empty throughout the entire paper.
	
	.
	\begin{definition}\label{d1}
		\cite{Chao}A real-valued function  $\hslash$ on $\aleph$ with(associated with) respect to(W.R.T. in short) the map $ b : \aleph \times \aleph \times [0, 1] \rightarrow \mathbf{R}$ is known to be {\it \bf sub-$b$-convex} if
		
		$$ \hslash(\sigma b_1 + (1- \sigma) b_2) \leq \sigma \hslash(b_1) + (1- \sigma) \hslash(b_2)+b(b_1,b_2, \sigma)$$ 
		holds $\forall ~~b_1, b_2 \in P $ and $\sigma \in [0, 1].$
	\end{definition}
	
	\begin{definition}\label{d2}
		\cite{Hudzik}A real-valued function  $\hslash$ on $\aleph$ is known to be {\it \bf $s$-convex in the second sense} if
		
		$$ \hslash(\sigma b_1 + (1- \sigma) b_2) < \sigma^s \hslash(b_1) + (1- \sigma)^s \hslash(b_2)$$ 
		holds $\forall ~~b_1, b_2 \in P $ , $\sigma \in [0, 1]$ and for some fixed $s \in (0,1].$
	\end{definition}
	\begin{definition}\label{d3}
		\cite{Liao}A real-valued function  $\hslash$ on $\aleph$ W.R.T. the map $ b : \aleph \times \aleph \times [0, 1] \rightarrow \mathbf{R}$ is known to be {\it \bf sub-$b$-$s$-convex in the second sense} if

		$$ \hslash(\sigma b_1 + (1- \sigma) b_2) < \sigma^s \hslash(b_1) + (1- \sigma)^s \hslash(b_2)+b(b_1,b_2, \sigma)$$ 
		holds $\forall ~~b_1, b_2 \in P, $  $\sigma \in [0, 1]$ and for some fixed $s \in (0,1].$
	\end{definition}
	Motivated by Definition \ref{d1}, \ref{d2} and \ref{d3}, we gives the new concepts of general $s$-convex function which is given below.
	\begin{definition}\label{d4}
		A real-valued function  $\hslash$ on $\aleph$ W.R.T. the function $ \vartheta :  \aleph \times [0, 1] \rightarrow \mathbf{R}$ is known to be {\it \bf general $s$-convex} if
		
		\begin{equation}\label{i1}
			\hslash(\sigma b_1 + (1- \sigma) b_2) \leq \sigma^s [\hslash(b_1)+\vartheta(b_1, \sigma)] + (1- \sigma)^s [\hslash(b_2)+\vartheta(b_2,\sigma)]+\vartheta(\frac{b_1+b_2}{2}, \sigma)	
		\end{equation}
		
		holds $\forall ~~b_1, b_2 \in P, $  $\sigma \in [0, 1]$ and for some fixed $s \in (0,1].$
	\end{definition}
	\begin{remark}\label{r1}
		If we take $\vartheta(b_1, \sigma)=0$, then general $s$-convex function  becomes $s$-convex in the second sense. On the oterhand, if we take $\vartheta(b_1, \sigma)=0$ and $s=1$, then general $s$-convex function  becomes convex.
	\end{remark}

	\begin{theorem}\label{t1}
		Suppose that $ \hslash_1, \hslash_2 : \aleph \rightarrow \mathbf{R}$ are general $s$-convex function W.R.T. the same map $ \vartheta : \aleph \times [0, 1] \rightarrow \mathbf{R}$, then $ \hslash_1+\hslash_2$ and $ \alpha \hslash_1 (\alpha \geq0)$ are general $s$-convex function W.R.T. the same map.
	\end{theorem}

	\begin{corollary}\label{c1}
		Suppose that $ \hslash_k : \aleph \rightarrow \mathbf{R}$, $(k=1,2,....,n)$ are general $s$-convex function W.R.T. the same map $ \vartheta_k : \aleph \times [0, 1] \rightarrow \mathbf{R}$, $(k=1,2,....,n)$ , respectively, then the function $$ \hslash=\sum_{k=1}^{n} \alpha_i\hslash_k, (\alpha_k \geq 0), (k=1,2,....,n)$$ is general $s$-convex function W.R.T. the  map  $ \vartheta=\sum_{k=1}^{n} \alpha_i\vartheta_k, , (k=1,2,....,n)$.
	\end{corollary}
	\begin{proposition}\label{p1}
		Suppose that $ \hslash_k : \aleph \rightarrow \mathbf{R}$, $(k=1,2,3,....,n)$ are general $s$-convex function W.R.T. the same map $ \vartheta_k : \aleph \times [0, 1] \rightarrow \mathbf{R}$, $(k=1,2,....,n)$ , respectively, then the function $ \hslash= \max\{\hslash_i, k=1,2,....,n\}$ is general $s$-convex function W.R.T. the  map  $ \vartheta=\max\{\vartheta_k, k=1,2,....,n\}$.
	\end{proposition}
	\begin{theorem}\label{t2}
		Suppose that $ \hslash_1 : \aleph\rightarrow \mathbf{R}$, are general $s$-convex function W.R.T. the map $ \vartheta : \aleph \times [0, 1] \rightarrow \mathbf{R}$ and the another function  $ \hslash_2 :  \mathbf{R} \rightarrow \mathbf{R}$ which is linear as well as non-decreasing, then $ \hslash_2 \circ \hslash_1$ is a general $s$-convex function W.R.T. the map $\hslash_2 \circ \vartheta.$
	\end{theorem}

	Now, we introduced a new theory of general $s$-convex set as below.
	\begin{definition}\label{t5}
		Assume $\aleph_1$ be a non-empty subset of $\mathbf{R^{m+1}}$. Then, the set $\aleph$ is known to be {\it \bf general $s$-convex set} W.R.T. the map $ \vartheta : \mathbf{R^m} \times [0, 1] \rightarrow \mathbf{R}$, if $$	(\sigma b_1 + (1- \sigma) b_2, ~\sigma^s [\alpha +\vartheta(b_1, \sigma )] + (1- \sigma)^s [\beta+\vartheta(b_2, \sigma)]+\vartheta(\frac{b_1+b_2}{2}, \sigma)) \in \aleph$$
		holds $\forall ~~(b_1, \alpha), (b_2, \beta ) \in \aleph_1, b_1, b_2 \in \mathbf{R^m}, \sigma \in [0, 1],$ and for some fixed $s \in (0,1]$.
	\end{definition}
	Here, we characterise the general $s$-convex function in terms of their epigraph $E(\aleph)$, which is provided by
	$$E(\hslash)=\{(b,\beta): b\in \aleph, \beta \in \mathbf{R}, \hslash(b)\leq \beta\}$$
	\begin{theorem}\label{t3}
		A function $ \hslash : \aleph \rightarrow \mathbf{R}$ is a general $s$-convex W.R.T. the map $\vartheta : \aleph \times [0, 1] \rightarrow\mathbf{R}$, $\iff$ $E(\hslash)$ is a general $s$-conevx set W.R.T. to $\vartheta$. 
	\end{theorem}

	\begin{proposition}\label{p2}
		Let us suppose that $\aleph_i$ is a family of general $s$-convex set W.R.T. the same map $\vartheta$, then $\cap_{i \in I}\aleph_i$ is also a general $s$-convex set W.R.T. the map $\vartheta$.
	\end{proposition}
	\begin{proof}
		Take $(b_1, \beta_1), (b_2, \beta_2) \in \cap_{i \in I}\aleph_i$, then, for any $i \in I, (b_1, \beta_1), (b_2, \beta_2) \in \aleph_i$. As $\aleph_i$ is a general $s$-convex set W.R.T. map $\vartheta$, for some $s \in (0, 1]$ and $\forall \sigma \in [0,1]$, it demonstrate that
		$$	(\sigma b_1 + (1- \sigma) b_2, ~\sigma^s [\beta_1 +\vartheta(b_1, \sigma )] + (1- \sigma)^s [\beta_2+\vartheta(b_2, \sigma)]+\vartheta(\frac{b_1+b_2}{2}, \sigma)) \in \aleph_i, $$ $\mbox{for any}~~ i \in I.$ Hence,  $$	(\sigma b_1 + (1- \sigma) b_2, ~\sigma^s [\beta_1 +\vartheta(b_1, \sigma )] + (1- \sigma)^s [\beta_2+\vartheta(b_2, \sigma)]+\vartheta(\frac{b_1+b_2}{2}, \sigma)) \in \cap_{i \in I}\aleph_i.$$ Therefore, $\cap_{i \in I}\aleph_i$ is a general $s$-convex set W.R.T. map $\vartheta$ and the proof is completed.
	\end{proof}
	\begin{proposition}\label{p3}
		If   $ \{ \hslash_i: i \in K\}$is a collection of arbitrary functions and all $\hslash_i$ is a general $s$-convex function W.R.T. the common map $\vartheta$, then the arbitrary function $\hslash=\sup_{i \in K}\hslash_i(b_1)$ is a general $s$-convex function W.R.T. common map $\vartheta$. 
	\end{proposition}

	\section{Main Results} Inside this section, we assume that $\hslash$ is a differentiable and general $s$-convex functions W.R.T. the map $\vartheta$. Also, we take the map $\vartheta(b, \sigma)$ such that $\dfrac{\vartheta(b, \sigma)}{\sigma}$ is well-defined.\\ Moreover, we suppose that limit $\displaystyle \lim_{\sigma \to  0^+}$ $\dfrac{\vartheta(\frac{b_1+b_2}{2}, \sigma)}{\sigma}$ exists and the limit is the maximum of $\dfrac{\vartheta\left(\frac{b_1+b_2}{2}, \sigma\right)-\vartheta(b_2, \sigma)-o(h)}{\sigma}$ $\forall \sigma \in (0, 1]$ and fixed $b_1,b_2 \in \aleph$ 
	\begin{theorem}\label{t4}
		Assume that a non-negative differential function $ \hslash : \aleph \rightarrow \mathbf{R}$ which is general $s$-convex W.R.T. the map $\vartheta$. Then\\
		
		\begin{eqnarray*}
			(a) \nabla \hslash(b_2)^T(b_1-b_2)&\leq& \sigma^{s-1}[\hslash(b_1)+\hslash(b_2)+\vartheta(b_1, \sigma)+\vartheta(b_2, \sigma)]+\displaystyle \lim_{\sigma \to  0^+} \dfrac{\vartheta(\frac{b_1+b_2}{2}, \sigma)}{\sigma}\\ 
			(b) \nabla \hslash(b_2)^T(b_1-b_2)&\leq& \sigma^{s-1}[\hslash(b_1)-\hslash(b_2)+\vartheta(b_1, \sigma)-\vartheta(b_2, \sigma)]\\&&+\dfrac{\hslash(b_2)}{\sigma}+\dfrac{\vartheta(b_2, \sigma)}{\sigma}+\displaystyle \lim_{\sigma \to  0^+} \dfrac{\vartheta\left(\frac{b_1+b_2}{2}, \sigma\right)}{\sigma} 
		\end{eqnarray*}
	\end{theorem}

	\begin{theorem}\label{t5}
		Assume that a non-negative differential function $ \hslash : \aleph \rightarrow \mathbf{R}$ which is general $s$-convex W.R.T. the non-postive map $\vartheta$. Then
		\begin{eqnarray*}
			\nabla \hslash(b_2)^T(b_1-b_2)&\leq& \sigma^{s-1}[\hslash(b_1)-\hslash(b_2)+\vartheta(b_1, \sigma)-\vartheta(b_2, \sigma)]+\displaystyle \lim_{\sigma \to  0^+} \dfrac{\vartheta(\frac{b_1+b_2}{2}, \sigma)}{\sigma}\\ 
		\end{eqnarray*}
	\end{theorem}

	\begin{corollary}\label{c2}
		Assume that a non-negative differential function $ \hslash : \aleph \rightarrow \mathbf{R}$ which is general $s$-convex W.R.T. the map $\vartheta$. Then
		\begin{eqnarray}\label{i16}
			\nabla(\hslash(b_2)-\hslash(b_1))^T(b_1-b_2) &\leq&  \dfrac{\hslash(b_1)}{\sigma}+\dfrac{\hslash(b_2)}{\sigma}+\dfrac{\vartheta(b_2,\sigma)}{\sigma}+\dfrac{\vartheta(b_2,\sigma)}{\sigma}\notag\\&&+2 \displaystyle \lim_{\sigma \to  0^+} \dfrac{\vartheta(\frac{b_1+b_2}{2}, \sigma)}{\sigma}.
		\end{eqnarray}
		If $\hslash$ is negative fucntion and $\vartheta$ is non-postive map, then
		\begin{eqnarray}\label{i17}
			\nabla(\hslash(b_2)-\hslash(b_1))^T(b_1-b_2) 
			&\leq&  2\displaystyle \lim_{\sigma \to  0^+}\dfrac{\vartheta(\frac{b_1+b_2}{2}, \sigma)}{\sigma}.
		\end{eqnarray}
	\end{corollary}

	We now apply the related findings from above to the nonlinear programming. First, we think about the unconstraint problem (S).
	\begin{equation}\label{i22}
		(S):	\min \{\hslash(b), b \in \aleph\}
	\end{equation}
	\begin{theorem}\label{t6}
		Suppose that $\hslash$ a real-valued non-negative differentiable function on $\aleph$ and general $s$-convex function associated with the map $\vartheta$. If $ b_2 \in \aleph$ and the inquality
		\begin{eqnarray}\label{i23}
			\nabla \hslash(b_2)^T(b_1-b_2)-\dfrac{\hslash(b_2)}{\sigma}-\dfrac{\vartheta(b_2, \sigma)}{\sigma}-\displaystyle \lim_{\sigma \to  0^+} \dfrac{\vartheta\left(\frac{b_1+b_2}{2}, \sigma\right)}{\sigma}&\geq& \sigma^{s-1}[\vartheta(b_1, \sigma)-\vartheta(b_2, \sigma)] \notag\\
		\end{eqnarray}
		holds $\forall b_1 \in \aleph, \sigma \in [0,1]$ and for some fixed $s \in (0,1]$, so $b_2$ is the optimal solution to the optimal problem (S) associated with $\hslash$ on $\aleph$.
	\end{theorem}

	Now let's look at the next general $s$-convex programming unconstraints
	\begin{equation}\label{i100}
		S:	\min \{\hslash(b), b \in [1, \infty)\}
	\end{equation}
	where  $\hslash(b)=[(b-1)^2+(b-1)]^s,$ where $s$ is a fixed number in $(0, 1)$,and  $\vartheta(b, \sigma)=\sigma(2b+6)$. As $\hslash$ be the non-negative differentiable and general $s$-convex function associated with the map $\vartheta$, and the limit $ \displaystyle \lim_{\sigma \to  0^+}\dfrac{\vartheta\left(\frac{b_1+b_2}{2}, \sigma\right)}{\sigma}$ exists for fixed $b_1, b_2 \in [1, \infty)$ and $\sigma \in (0,1]$. Calculating after that, here's what we have
	$$ \nabla\hslash(b^*)^T(b_1-b^*)=s[(b^*-1)^2+(b^*-1)]^{s-1}(2(b^*-1)+1)(b_1-b^*)$$
	$$\dfrac{\hslash(b_1)}{\sigma}=\dfrac{[(b-1)^2+(b-1)]^s}{\sigma}$$
	$$\displaystyle \lim_{\sigma \to  0^+}\dfrac{\vartheta\left(\frac{b_1+b_2}{2}, \sigma\right)}{\sigma}=b_1-b^*$$
	It is easy to see that at $b^*=1$, the inequality
	\begin{eqnarray*}
		\nabla \hslash(b_2)^T(b_1-b_2)-\dfrac{\hslash(b_2)}{\sigma}-\dfrac{\vartheta(b_2, \sigma)}{\sigma}-\displaystyle \lim_{\sigma \to  0^+} \dfrac{\vartheta\left(\frac{b_1+b_2}{2}, \sigma\right)}{\sigma}&\geq& \sigma^{s-1}[\vartheta(b_1, \sigma)-\vartheta(b_2, \sigma)] \notag\\
	\end{eqnarray*}holds $\forall b_1 \in [1, \infty), \sigma \in (0, 1]$ and some fixed number $s \in (0, 1)$. So, by Theorem \ref{t6}, at $b_1=1$, $\hslash(b_1)$ gives the minimum value.
	\begin{corollary}\label{c3}
		Suppose that $\hslash$ a real-valued non-negative differentiable function on $\aleph$ and strictly general $s$-convex function associated with the map $\vartheta$. If the condition \ref{i23} satisfies and $b_2 \in \aleph$, then $b_2$ is optimal solution which is unique.
	\end{corollary}

	The resulting results are then applied in the following ways to nonlinear programming with inequality constraints:
	
	\begin{equation}\label{i26}
		(S_p):	\min \{\hslash(b)| b\in \mathbf{R^m},~ b \in \aleph, f_i(b)\leq0, i \in I\}, I=\{1,2,3,....,n\}
	\end{equation}
	Let $E= \{ b\in \mathbf{R^m}: f_i(b)\leq0, i \in I\}$. We assume for the sake of explanation that $\hslash$ and $f_i$ are both differentiable and that $E$ is a nonempty set in $\mathbf{R^m}$.
	
	\begin{theorem}\label{t7}
		\{ Karush-Kuhn-Tucker Sufficint Conditions\} Let $\hslash : \mathbf{R^m} \rightarrow \mathbf{R}$ be a non-negative differentiable general $s$-convex function associated with the map $\vartheta : \mathbf{R^m} \times (0,1] \rightarrow \mathbf{R}$ and $f_i : \mathbf{R^m} \rightarrow \mathbf{R}$ $ (i \in I)$ are differentiable general $s$-convex function associated with the map $\vartheta_i : \mathbf{R^m} \times (0,1] \rightarrow \mathbf{R}$ $(i \in I)$. Let $b^* \in E$ is a KKT point of $(S_p),$i.e., $\exists$ multipliers $v_i\geq0$ $(i \in I)$ s.t.
		\begin{eqnarray}\label{i27}
			\nabla\hslash(b^*)+ \sum_{i \in I}v_i\nabla f_i(b^*)=0,~~ v_if_i(b^*)=0.
		\end{eqnarray}
		Assume that the inequality
		\begin{eqnarray}\label{i28}
			\dfrac{\phi(b^*)}{\sigma}+\dfrac{\vartheta(b^*, \sigma)}{\sigma}+\displaystyle \lim_{\sigma \to  0^+} \dfrac{\vartheta\left(\frac{t_1+b^*}{2}, \sigma\right)}{\sigma}&\leq&-\sum_{i \in I}v_i\displaystyle \lim_{\sigma \to  0^+} \dfrac{\vartheta\left(\frac{t_1+b^*}{2}, \sigma\right)}{\sigma}\notag\\&&-2\sigma^{s-1}[\vartheta(t_1, \sigma)-\vartheta(b^*, \sigma)] ,
		\end{eqnarray}
		holds $\mbox{for all}~~ b, b^* \in \mathcal{R}^m$, then $b^*$ be an optimal solution of the problem $(S_p)$.
	\end{theorem}

\end{document}